# Identities Inspired by the Ramanujan Notebooks

# Second Series

by Simon Plouffe
First draft August 2006
Revised March 14, 2011

## Abstract

A series of formula is presented that are all inspired by the Ramanujan Notebooks [6]. One of them appears in the notebooks II

$$\zeta(3) = \frac{7\pi^3}{180} - 2\sum_{n=1}^{\infty}\frac{1}{n^3(e^{2\pi n}-1)}$$

That formula inspired others that appeared in 1998, 2006 and 2009 on the author's website and later in literature [1][2][3]. New formulas for $\pi$ and the Catalan constant are presented along with a surprising series of approximations. A new set of identities is given for Eisenstein series. All of the formulas are conjectural since they were found experimentally. A new method is presented for the computation of the partition function.

Une série de formules utilisant l'exponentielle est présentée. Ces résultats reprennent ceux apparaissant en 1998, 2006 et 2009 sur [1][2][3]. Elles sont toutes inspirées des Notebooks de Ramanujan tels que

$$\zeta(3) = \frac{7\pi^3}{180} - 2\sum_{n=1}^{\infty}\frac{1}{n^3(e^{2\pi n}-1)}$$

Une nouvelle série pour $\pi$ et la constante de Catalan sont présentés ainsi qu'une série d'approximations surprenantes. Une série d'identités nouvelles sont présentées concernant les séries d'Eisenstein. Toutes les formules présentées sont des conjectures, elles ont toutes été trouvées expérimentalement. Une nouvelle méthode est présentée pour le calcul des partages d'un entier.

# 1. Introduction

By taking back the series found in 2006, I extended the search to more general expressions with exponents 1,2 and 4 for the exponential term and found the following

The same pattern is present for powers of $\pi$ and $\zeta(n)$

1.1 $$\pi = 72 \sum_{n=1}^{\infty} \frac{1}{n(e^{\pi n} - 1)} - 96 \sum_{n=1}^{\infty} \frac{1}{n(e^{2\pi n} - 1)} + 24 \sum_{n=1}^{\infty} \frac{1}{n(e^{4\pi n} - 1)}$$

1.2 $$\frac{1}{\pi} = 8 \sum_{n=1}^{\infty} \frac{n}{e^{\pi n} - 1} - 40 \sum_{n=1}^{\infty} \frac{n}{e^{2\pi n} - 1} + 32 \sum_{n=1}^{\infty} \frac{n}{e^{4\pi n} - 1}$$

1.3 $$\pi^3 = 720 \sum_{n=1}^{\infty} \frac{1}{n^3(e^{\pi n} - 1)} - 900 \sum_{n=1}^{\infty} \frac{1}{n^3(e^{2\pi n} - 1)} + 180 \sum_{n=1}^{\infty} \frac{1}{n^3(e^{4\pi n} - 1)}$$

1.4 $$\zeta(3) = 28 \sum_{n=1}^{\infty} \frac{1}{n^3(e^{\pi n} - 1)} - 37 \sum_{n=1}^{\infty} \frac{1}{n^3(e^{2\pi n} - 1)} + 7 \sum_{n=1}^{\infty} \frac{1}{n^3(e^{4\pi n} - 1)}$$

1.5 $$\pi^5 = 7056 \sum_{n=1}^{\infty} \frac{1}{n^5(e^{\pi n} - 1)} - 6993 \sum_{n=1}^{\infty} \frac{1}{n^5(e^{2\pi n} - 1)} + 63 \sum_{n=1}^{\infty} \frac{1}{n^5(e^{4\pi n} - 1)}$$

1.6 $$\zeta(5) = 24 \sum_{n=1}^{\infty} \frac{1}{n^5(e^{\pi n} - 1)} - \frac{259}{10} \sum_{n=1}^{\infty} \frac{1}{n^5(e^{2\pi n} - 1)} - \frac{1}{10} \sum_{n=1}^{\infty} \frac{1}{n^5(e^{4\pi n} - 1)}$$

1.7 $$\pi^7 = \frac{907200}{13} \sum_{n=1}^{\infty} \frac{1}{n^7(e^{\pi n} - 1)} - 70875 \sum_{n=1}^{\infty} \frac{1}{n^7(e^{2\pi n} - 1)} + \frac{14175}{13} \sum_{n=1}^{\infty} \frac{1}{n^7(e^{4\pi n} - 1)}$$

1.8 $$\zeta(7) = \frac{304}{13} \sum_{n=1}^{\infty} \frac{1}{n^7(e^{\pi n} - 1)} - \frac{103}{4} \sum_{n=1}^{\infty} \frac{1}{n^7(e^{2\pi n} - 1)} + \frac{19}{52} \sum_{n=1}^{\infty} \frac{1}{n^7(e^{4\pi n} - 1)}$$

For the Catalan constant, I find this new identity:

1.9 $$K = 11 \sum_{n=1}^{\infty} \frac{1}{n^2(cosh\ (\pi n) - 1)} - \frac{71}{2} \sum_{n=1}^{\infty} \frac{1}{n^2(cosh\ (2\pi n) - 1)} + 11 \sum_{n=1}^{\infty} \frac{1}{n^2(cosh\ (4\pi n) - 1)}$$

For $1/\pi^2$, by varying the function at the numerator, I find this:

1.10
$$\frac{1}{\pi^2} = 4 \sum_{n=1}^{\infty} \frac{\sigma_1(n)n}{e^{\pi n}} - 64 \sum_{n=1}^{\infty} \frac{\sigma_1(n)n}{e^{2\pi n}} + 64 \sum_{n=1}^{\infty} \frac{\sigma_1(n)n}{e^{4\pi n}}$$

Here, $\sigma_1(n)$ is Euler's sigma function of order 1. Actually the same coefficients as with $cosh\ (k\pi n)$, k=1,2 and 4.

1.11
$$\frac{1}{\pi^2} = 2 \sum_{n=1}^{\infty} \frac{n^2}{cosh(\pi n) - 1} - 32 \sum_{n=1}^{\infty} \frac{n^2}{cosh\ (2\pi n) - 1} + 32 \sum_{n=1}^{\infty} \frac{n^2}{cosh\ (4\pi n) - 1}$$

The pattern persist for $1/\pi^3$ but apparently for no other powers of $\pi$

1.12
$$\frac{1}{\pi^3} = 4 \sum_{n=1}^{\infty} \frac{\sigma_1(n)n^2}{e^{\pi n}} - 128 \sum_{n=1}^{\infty} \frac{\sigma_1(n)n^2}{e^{2\pi n}} + 256 \sum_{n=1}^{\infty} \frac{\sigma_1(n)n^2}{e^{4\pi n}}$$

## 2. Experiments with fractional exponent

I was compiling a table of values for the Inverter [9] and found that for some arguments the closeness to rational numbers. These are the 2 examples that are the most striking.

2.1
$$\sum_{n=1}^{\infty} \frac{n^3}{e^{2\pi n/7} - 1} = 10.00000000000000190161767888663\ ...$$

2.2
$$\sum_{n=1}^{\infty} \frac{n^3}{e^{2\pi n/13} - 1} \cong 119.000000000000000000000000000000959374585\ ...$$

The precision is 15 and 31 decimal digits for an argument of $2\pi n/163$ the precision is 435 decimal digits. Other series of the form $\sum_{n=1}^{\infty} \frac{n^3}{e^{\frac{2\pi n}{k}}-1}$ are also producing near integers when k is not a multiple of 2,3 and 5. For the exponent one can obtain near integers when the exponent of n is 4m-1, m > 0. This fact is related to properties of Eisenstein series which is; if $240\ |\ k^4 - 1$ then the series is near an integer. But it does not always produces approximations since I have this identity for $\pi$. We see the pattern [1,2,4] again with the exponent.

2.3
$$\frac{\pi}{10} = -\sum_{n=1}^{\infty} \frac{1}{n(e^{\pi n} - 1)} + 4 \sum_{n=1}^{\infty} \frac{1}{n(e^{2\pi n} - 1)}$$
$$-1 \sum_{n=1}^{\infty} \frac{1}{n(e^{4\pi n} - 1)} + \sum_{n=1}^{\infty} \frac{1}{n(e^{\pi n/5} - 1)} - 4 \sum_{n=1}^{\infty} \frac{1}{n(e^{2\pi n/5} - 1)}$$
$$+ \sum_{n=1}^{\infty} \frac{1}{n(e^{4\pi n/5} - 1)}$$

2.4
$$\frac{7\pi}{120} = -2 \sum_{n=1}^{\infty} \frac{1}{n(e^{\pi n} - 1)} - \sum_{n=1}^{\infty} \frac{1}{n(e^{\pi n/5} - 1)} + 4 \sum_{n=1}^{\infty} \frac{1}{n(e^{2\pi n/5} - 1)} - 1 \sum_{n=1}^{\infty} \frac{1}{n(e^{4\pi n/5} - 1)}$$

$$2.5 \quad 3\log(\varphi) = -4\sum_{n=1}^{\infty}\frac{1}{n(e^{\pi n}-1)} + 10\sum_{n=1}^{\infty}\frac{1}{n(e^{2\pi n}-1)}$$
$$-4\sum_{n=1}^{\infty}\frac{1}{n(e^{4\pi n}-1)} + 4\sum_{n=1}^{\infty}\frac{1}{n(e^{\frac{\pi n}{5}}-1)} - 10\sum_{n=1}^{\infty}\frac{1}{n(e^{2\pi n/5}-1)}$$
$$+ 4\sum_{n=1}^{\infty}\frac{1}{n(e^{4\pi n/5}-1)}$$

$$2.6 \quad \sum_{n=1}^{\infty}\frac{1}{n(e^{\pi n/5}-1)} - 4\sum_{n=1}^{\infty}\frac{1}{n(e^{2\pi n/5}-1)} + \sum_{n=1}^{\infty}\frac{1}{n(e^{4\pi n/5}-1)}$$
$$= \frac{\pi}{40} - \frac{3\ln(\pi)}{2} + 2\ln\Gamma\left(\frac{1}{4}\right) - \frac{7\ln(2)}{4}$$

$$2.7 \quad \log(\varphi) = \sum_{n=1}^{\infty}\frac{1}{n(e^{\frac{\pi n}{5}}-1)} - 2\sum_{n=1}^{\infty}\frac{1}{n(e^{\frac{2\pi n}{5}}-1)} + \sum_{n=1}^{\infty}\frac{1}{n(e^{\frac{4\pi n}{5}}-1)} + \frac{\pi}{120} - \frac{\ln(2)}{4}$$

Where $\varphi$ is the golden ratio.

Since the generic series for $\pi$ is with $\sum_{n=1}^{\infty}\frac{1}{n(e^{\pi n}-1)}$ that series is the log of the well known Euler partition function $F(x) = \prod_{n=1}^{\infty}\frac{1}{1-x^n}$ series which means that an identity can be translated into the partition function as well, of course when $x \to e^{-\pi x}$.

$$2.8 \quad \frac{F(1/5)^5 F(4/5)^5}{F(2/5)^{20}} = \frac{F(2)^{28}}{F(1)^{31}F(4)^7}$$

The exact expression for $\sum_{n=1}^{\infty}\frac{1}{n(e^{\pi n}-1)}$ can also be found by a variety of methods. One of them is to simply try it in the Integer Relations algorithm to find (2.6). Since the partition function is involved it is natural to ask if there could be a relationship with the Roges-Ramanujan, recall G(x), H(x) are and let's define $J(x)$,

$$2.9 \quad J(x) = G(x)H(x) = \prod_{n=1}^{\infty}(1-x^{(5n-1)})(1-x^{(5n-4)})(1-x^{(5n-2)})(1-x^{(5n-3)})$$

And again with $x \to e^{-\pi x}$, then

$$2.10 \quad \frac{J(1/5)^{10}J(4/5)^{10}}{J(2/5)^{40}} = e^{\pi}$$

$$2.11 \quad J(2/5)^{30} = \frac{e^{2\pi}}{\varphi^{15}}$$

$$2.12 \quad \frac{J(1)^6 J(4)^{10}}{J(2)^{12}} = \frac{e^{\pi}}{\varphi^6}$$

$$2.13 \quad \frac{J(1/5)^4 J(4/5)^4}{J(2/5)^{10}} = \varphi^3$$

$$2.14 \quad J(1/5)^6 J(4/5)^6 = e^{\pi}\varphi^{12}$$

The exact expression for k=2/7 was found by Bill Gosper using the Computer Algebra system Macsyma and a set of personal routines.

$$2.15 \qquad \sum_{n=1}^{\infty} \frac{n^3}{e^{2\pi n/7} - 1} = \frac{-1}{240} + \frac{1}{320}(301 + 210\sqrt{2}\, 7^{1/4} + 120\sqrt{7} + 90\sqrt{2}\, 7^{3/4}) \frac{\pi^2}{\Gamma(\frac{3}{4})^8}$$

In fact that series is the series of Eisenstein which are.

$$2.16 \qquad E_4(q) = 1 + 240 \sum_{k=1}^{\infty} \sigma_3(k) q^{2k}$$

$$2.17 \qquad E_8(q) = 1 + 480 \sum_{k=1}^{\infty} \sigma_7(k) q^{2k}$$

$$2.18 \qquad E_{12}(q) = 1 + \frac{65520}{691} \sum_{k=1}^{\infty} \sigma_{11}(k) q^{2k}$$

When with $q \to e^{-\pi q}$, $E_4(1/10)$, $E_4(1/5)$, $E_4(2/5)$ we can get respectively

$$2.19 \qquad \frac{\pi^2}{\Gamma(3/4)^8}\left(\frac{5313}{4} + 630\sqrt{5} + 90\sqrt{360 + 161\sqrt{5}}\right)$$

$$2.20 \qquad \frac{\pi^2}{\Gamma(3/4)^8}\left(\frac{(483}{4} 90\sqrt{5})\right)$$

$$2.21 \qquad \frac{\pi^2}{\Gamma(3/4)^8}\left((5313 + 2520\sqrt{5} - 32\sqrt{\frac{91125}{2} + \frac{326025\sqrt{5}}{16}})\right)$$

Note, $E_4(1/10) = 10000.000000000000000000000012378...$

$$2.22 \qquad 2\sum_{n=1}^{\infty} \frac{n^3}{e^{\frac{\pi n}{5}} - 1} - 28 \sum_{n=1}^{\infty} \frac{n^3}{e^{\frac{2\pi n}{5}} - 1} + 32 \sum_{n=1}^{\infty} \frac{n^3}{e^{\frac{4\pi n}{5}} - 1} + 28 \sum_{n=1}^{\infty} \frac{n^3}{e^{\frac{2\pi n}{5}} - 1} - 257 \sum_{n=1}^{\infty} \frac{n^3}{e^{\pi n} - 1} + 251 \sum_{n=1}^{\infty} \frac{n^3}{e^{2\pi n} - 1} = 0$$

$$2.23 \qquad 8 \sum_{n=1}^{\infty} \frac{n^7}{e^{\pi n/5} - 1} - 2192 \sum_{n=1}^{\infty} \frac{n^7}{e^{2\pi n/5} - 1} + 2048 \sum_{n=1}^{\infty} \frac{n^7}{e^{4\pi n/5} - 1} + 208897 \sum_{n=1}^{\infty} \frac{n^7}{e^{\pi n} - 1} - 208761 \sum_{n=1}^{\infty} \frac{n^7}{e^{2\pi n} - 1} = 0$$

The pattern observed earlier [1,2,4] can be translated into Eisenstein series identities that are new (?). Here $E_n(q)$ is with $q \to e^{-2\pi q}$

2.24 $\quad -E_4(1/10)+14E_4(1/5)-16E_4(2/5)+\dfrac{1288}{11}E_4(1/2)=0$

2.25 $\quad E_8(1/10)-274E_8(1/5)+256E_8(2/5)+\dfrac{3133472}{121}E_8(1/2)=0$

2.26 $\quad -E_{12}(1/10)+4034E_{12}(1/5)-4096E_{12}(2/5)+\dfrac{7811747968}{2081}E_{12}(1/2)=0$

I could not find any other exact fractional values, the other ones are only translatable into approximations like, here $F$ as in (2.8). Other values are of interest as well like

|  | $F(x)=\displaystyle\prod_{n=1}^{\infty}\dfrac{1}{1-x^n}\,,\,x\to e^{-\pi x}$ | Precision in digits |
|---|---|---|
| 2.27 | $F(1)^8\approx\dfrac{e^{\pi}}{2^4}$ | 4 |
| 2.28 | $F(1/2)^{16}\approx\dfrac{e^{5\pi}}{2^{16}}$ | 6 |
| 2.29 | $\dfrac{F(1/8)^{36}F(1/15)^{18}}{F(1/12)^{36}F(3/20)^{18}}\approx e^{\pi}$ | 35 |
| 2.30 | $F(1/8)^{64}\approx\dfrac{e^{85\pi}}{2^{128}}$ | 36 |
| 2.31 | $\dfrac{F(1/12)^4F(1/36)^2}{F(1/9)^2F(1/24)^4}\approx e^{\pi}$ | 48 |
| 2.32 | $F(1/32)^{256}\approx\dfrac{e^{1365\pi}}{2^{768}}$ | 173 |

# 3. The computation of $p(n)$

A brief history note : In 1917, a certain Percy Alendander MacMahon computed by hand the value of p(200) in about six months. The value is 3972999029388. Later, G.H. Hardy and S. Ramanujan used the so-called circle method to compute p(200) using few terms of an asymptotic formula :

3.1
$$p(n) \approx A_n e^{\pi\sqrt{\frac{2}{3}(n-\frac{1}{24})}}$$

Where $A_n = \frac{1}{2n\sqrt{2}}\left(\frac{\pi}{\sqrt{6(\frac{n-1}{24})}} - \frac{1}{2(\frac{n-1}{24})^{3/2}}\right)$. The later improvement of Rademacher is the one used in the Computer Algebra System Mathematica for small values of n.

One of the few experimental data I have found are dealing with the partition function. For example, the number $F(4) = \frac{2^{\frac{7}{8}}\Gamma(\frac{3}{4})}{\pi^{1/4}e^{\pi/6}} = 1.00000348736667944964958540340...$ and the number

3.2
$$F(10) = \frac{\Gamma(\frac{3}{4})\sqrt{5+5\sqrt{5}}}{\pi^{1/4}e^{\pi/6}} = 1.00000000000002271101068...$$

Can easily be found using Integer Relations Algorithms.

In other words, since the evaluation of $F(x) = \prod_{n=1}^{\infty}\frac{1}{1-x^n}$ is made at $x \to e^{-\pi x}$, then it also means that the number $= 1.00000000000002271101068...$ when expanded into the base $e^{\pi}$ would normally produce the numbers $p(n)$, given enough numerical precision of course.

That idea of using a certain base to expand numbers if not new, I used it in a paper of 1993 to get the algebraic generating function of a set of sequences [12]. The idea is based on one observation: The number $\frac{\sqrt{51}}{14} - \frac{1}{2} = 0.0101020306102035...$ and the coefficients are 1, 1, 2, 3, 6, 10, 20, 35, ... are related to the central binomial coefficients.

This is not surprising since for example, the function $\frac{1}{\sqrt{1-4x}}$ when evaluated at x=1/100 is equal to 1.02062072615965754091.... and $\frac{1}{\sqrt{1-4x}}$ is the algebraic generating function of the binomial coefficients, we can *see* the coefficients.

The question that comes naturally after is what does happen we use the base $e^{\pi}$ for certain functions? The answer is the same thing happens but we need more numerical precision and we would not be able to see the coefficients as easily.

By using 2800 decimal digits, the number $F(10) = \frac{\Gamma(\frac{3}{4})\sqrt{5+5\sqrt{5}}}{\pi^{1/4} e^{\pi/6}}$ expanded in base $e^\pi$ will produce the first 205 terms of $F(x)$.

In particular $p(200) = 3972999029388$. To expand a number in base k we apply the formula :

3.3 $$y_n = [kx_n], x_{n+1} = \{kx_n\}$$

Where k is the base, $x_0$ is the real to be expanded, [ ] is the integer part function and { } is the fractional part of a real number.

Note: The computation of $p(200)$ in 1917 by the Major Percy MacMahon took 6 months, the computation was done by hand.

## 4. Conclusion

As far as the author knows the formulas with arguments [1/5, 2/5, 4/5] as well as the formula for the Catalan constant , π and 1/ π are new. No other exact formula was found for the fractional exponent but the search was limited to the Farey set of order 60. The approximation (2.24) for $e^\pi$ is the simplest found.